\documentclass[11pt]{amsart}
\usepackage{geometry}                
\usepackage{graphicx}
\usepackage{amssymb}
\usepackage{epstopdf}
\theoremstyle{plain}
\newtheorem{theo}{Theorem}[section]

\newtheorem{prop}[theo]{Proposition}

\theoremstyle{definition}

\def\D{\Delta}

\def\g{\gamma}

\def\l{\lambda}

\def \R{\mathbb R}

\def\W{\mathcal W}

\def\wt{\widetilde}
\def\({\biggl(}
\def\){\biggr)}

\def\<{\bold\langle}
\def\>{\bold\rangle}

\def\M{\widetilde {M}}
\def\MM{\widetilde {M}(\infty)}

\begin{document}
\title{Volume entropy rigidity of non-positively curved symmetric spaces}
\author{Fran\c cois Ledrappier}
\address{LPMA, UMR CNRS 7599, Universit\'e Paris 6, Bo\^ite Courrier 188, 4, Place Jussieu, 75252 PARIS cedex 05, France}
\address{Department of Mathematics, University of Notre Dame, Notre Dame, IN 46556, USA} \email{fledrapp@nd.edu}


\maketitle
\centerline{\it{To Werner Ballmann for his 60th birthday}}
\begin{abstract}
We characterize symmetric spaces of non-positive curvature by the equality case of general inequalities between geometric quantities.
\end{abstract}
\section{Introduction}

Let $(M,g)$ be a closed connected Riemannian manifold, and $\pi : (\M, \wt g) \to (M,g)$ its universal cover endowed with the lifted Riemannian metric. We denote $p(t,x,y), t \in \R_+, x,y \in \M $ the heat kernel on $\M$, the fundamental solution of the heat equation $ \frac{\partial u }{\partial t} = {\textrm { Div }}\nabla u $ on $\M$. Since we have a compact quotient, all the following limits exist as $t \to \infty $ and are  independent of $x \in \M$:
\begin{eqnarray*}
\l_0 \;&=&\;  \inf _{f \in C^2_c(\M)} \frac {\int |\nabla f|^2}{\int |f|^2} \; = \; \lim _t - \frac{1}{t} \ln p(t,x,x) \\
\ell \; &=&\;  \lim _t \frac{1}{t} \int d(x,y) p(t,x,y) d{\textrm {Vol}}(y) \\
h \; &=&\; \lim_t - \frac{1}{t} \int p(t,x,y) \ln p(t,x,y) d{\textrm {Vol}}(y) \\
v\; &=& \; \lim _t \frac{1}{t} \ln {\textrm {Vol}} B_{\M} (x,t) ,
\end{eqnarray*}
where $B_{\M} (x,t) $ is the ball of radius $t$ centered at $x$ in $\M$ and ${\textrm {Vol}}$ is the Riemannian volume on $\M$.

All these numbers are nonnegative. Recall $\l_0$ is the Rayleigh quotient of $\M$, $\ell $ the linear drift, $h$ the stochastic entropy and $v$ the volume entropy. There is the following relation:
\begin{equation}\label{ineq}
4\l_0 \; \stackrel {(a)}{\leq} \; h \;\stackrel {(b)}{ \leq }\; \ell v \;\stackrel{(c)}{ \leq }v^2.
\end{equation}
See \cite{L1} for (a), \cite{Gu} for (b). Inequality (c) is shown in \cite{L3} as a corollary of (b) and (\ref{basic}):
\begin{equation}\label{basic}
\ell ^2 \; \leq \; h 
\end{equation}
 If $(\M, g) $ is a locally symmetric space of nonpositive curvature, all five numbers $4\l_0, \ell ^2, h, \ell v$ and $ v^2$ coincide and are positive unless $(\M, g )$ is $(\R^n, {\textrm {Eucl.}})$. Our result is a partial converse:

\begin{theo}\label{main}
Assume $(M,g) $ has nonpositive curvature. With the above notation, any of the equalities
$$ \ell \; = \; v, \quad h \; = \; v^2 \; \; {\textrm { and }} \; \; 4\l_0 \; = \; v^2 $$ hold if, and only if, $(\M, \wt g)$ is a symmetric space.
\end{theo} 

As recalled in \cite{L3}, Theorem \ref{main} is known in negative curvature and follows from \cite{K}, \cite{BFL}, \cite{FL}, \cite{BCG} and \cite{L1}. The other possible converses are delicate:
even for negatively curved manifolds, in dimension greater than two, it is not known that $h = \ell v$ holds only for locally symmetric spaces. This is equivalent to a conjecture of Sullivan (see \cite{L2} for a discussion). Sullivan conjecture holds for surfaces of negative curvature (\cite {L1}, \cite {Ka}). It is not known either whether $4 \l_0 = h$ holds only for locally symmetric spaces. This would follow from the hypothetical $4 \l_0 \stackrel {(d)}{\leq }\ell ^2$ by  the arguments  of this note.

\

We assume henceforth that $(M,g)$ has nonpositive sectional curvature. Given a geodesic $\g $ in $M$, Jacobi fields along $\g$ are vector fields $ t \mapsto J(t) \in T_{\g(t)}M$ which describe infinitesimal variation of geodesics around $\g$. By nonpositive curvature, the function $t \mapsto \| J(t) \|$ is convex. Jacobi fields along $\g$ form a vector space of dimension $2$ Dim $M$. The rank of the geodesic $\g$ is the dimension of the space of Jacobi fields such that $t \mapsto \|J(t)\| $ is a constant function on $\R$. The rank of a geodesic $\g$ is at least one because of the trivial $t \mapsto \dot {\g}(t)$ which describes the variation by sliding the geodesic along itself.  The rank of the manifold $M$ is the smallest rank of geodesics in $M$. Using  rank rigidity theorem (\cite {B1}, \cite{BS}), we reduce in section 2 the proof of Theorem \ref{main} to proving that if $(M,g)$ is rank one, equality in (\ref{basic}) implies that $(\M, \wt g) $ is a symmetric space. For this, we show in section 3 that equality in (\ref{basic}) implies that $(\M, \wt g) $ is asymptotically harmonic (see the definition below). This uses the Dirichlet property at infinity (Ballmann \cite{B2}). Finally, it was recently observed by A. Zimmer (\cite{Z}) that asymptotically harmonic universal covers of rank one manifolds are indeed symmetric spaces.

\

\section{Generalities and reduction of Theorem \ref{main}}

\

We recall the notations and results from Ballmann's monograph \cite{B3} about the Hadamard manifold $(\M, \wt g)$ that we use. The space $\M$ is homeomorphic to a ball. The covering group $G := \pi _1(M) $ satisfies the duality condition (\cite {B3} page 45).

\subsection{Boundary at infinity}

Two geodesic rays $\g, \g'$ in $\M$ are said to be asymptotic if $\sup _{t \geq 0 } d(\g(t), \g'(t) ) < \infty $. The set of classes of asymptotic unit speed  geodesic rays is called the boundary at infinity $\MM$. $\M \cup \MM$ is endowed with the topology of a compact space where $\MM$ is a sphere and where, for each unit speed geodesic ray $\g$, $\g(t) \to [\g] $ as $t \to \infty $. The action of the group $G$ on $\MM$ is the continuous extension of its action on $\M$. For any $x, \xi \in \M \times \MM$, there is a unique unit speed geodesic $\g_{x, \xi}$ such that $\g_{x, \xi}(0) = x $ and $[\g_{x, \xi}] = \xi$. The mapping $\xi \mapsto \dot\g_{x, \xi}(0) $ is a homeomorphism $\pi_x ^{-1} $ between $\MM$ and the unit sphere $S_x\M$ in the tangent space at $x$ to $\M$. We will identify $S\M$ with $\M \times \MM$ by $(x,v) \mapsto (x, \pi _x v)$. Then the quotient $SM$ is identified with the quotient of $\M \times \MM$ under the diagonal action of $G$.

Fix $x_0 \in \M$ and $\xi  \in \MM $. The {\it  {Busemann function }} $b_\xi  $ is the function on $\M$ given by:
$$ b_\xi (x) \; = \; \lim _{y \to \xi } d(y,x) - d(y,x_0) .$$
Clearly, $b_{g\xi} (gx) = b_\xi (x) + b_{g\xi} (g x_0)$. Moreover, the function $x \mapsto b_\xi (x) $ is of class $C^2$ (\cite{HI}). It follows that the fonction  $\D _x b_\xi $ satisfies $\D _{gx}b_{g\xi} =  \D _xb_\xi $ and therefore defines a function $B$ on $ G \setminus (\M \times \MM ) = SM$. It follows from the argument of \cite{HI} that the function $B$ is continuous on $SM$ (see \cite{B3}, Proposition 2.8, page 69).

\subsection{Jacobi fields}

Let $(x,v) $ be  a point in $T\M$. Tangent vectors in $T_{x,v}T\M$ correspond to variations of geodesics and can be represented by Jacobi fields along the unique geodesic  $\g_{x,v}$ with initial value $\g(0) = x, \dot \g(0) = v$. A Jacobi field $J(t), t \in \R$ along $\g_{x,v}$ is uniquely determined by the values of $J(0)$  and  $J'(0)$. We describe tangent vectors in $T_{x,v}T\M$ by the associated pair $(J(0), J'(0))$ of vectors in $T_x\M$. The metric on $T_{x,v}T\M$ is given by $\|(J_0,J'_0)\| ^2 = \|J_0\| ^2 + \|J'_0\| ^2$. Assume $(x,v ) \in SM$. A vertical vector in $T_{x,v}S\M$ is a vector tangent to $S_x\M$.  It corresponds to a pair $(0, J'(0))$, with $J'(0) $ orthogonal to $v$. Horizontal vectors correspond to pairs $(J(0),0)$. In particular, let $X$ be the vector field on $S\M$ such that  the integral flow of $X$  is the geodesic flow. The geodesic spray $X_{x,v}$ is the horizontal vector associated to $(v,0)$. The orthogonal space to $X $ is preserved by the differential  $Dg_t $ of the geodesic flow. More generally, the Jacobi fields representation of $TT\M$ satisfies $D_{x,v} g_t (J(0),J'(0)) = (J(t), J'(t))$.

\

For any vector $Y \in T_x\M$, there is a unique vector $Z = S_{x,v} Y $ such that the Jacobi field $J$ with $J(0) = Y, J'(0) = Z $ satisfies $\|J(t) \| \leq C $ for $t \geq 0$ (\cite{B3} Proposition 2.8 (i)). The mapping $S_{x,v} : T_x\M \to T_x\M $ is linear and  selfadjoint. The vectors $(Y, SY)$ describe  variations of asymptotic geodesics and the  subspace $E^s_{x,v}  \subset T_{x,v}T\M $ they generate corresponds to $TW^s_{x,v}$, where $W^s_{x,v}$, the set of initial  vectors of geodesics asymptotic to $\g_{x,v}$, is identified with $\M \times \pi_x(v) $ in $\M \times \MM$.  Observe that $S_{x,v} \dot \g_{x,v}(0) = 0 $ and that the operator $S_{x,v}$ preserves $(\dot \g_{x,v}(0)) ^\perp .$ Recall from \cite{B3}, Proposition 3.2 page 71, that, for $Y \in (\dot\g_{x,v}(0))^\perp$, with $\pi _x v = \xi,$
$$ D_Y(\nabla  b_\xi ) \; = \; - S_{x,v} Y ,$$ and therefore $\D_x b_\xi  = - {\textrm { Tr }} S_{x,v}$ with $\pi _x(v) = \xi.$

Similarly, there is a selfadjoint linear operator $U_{x,v}  : T_x\M \to T_x\M $ such that  the Jacobi field $J$ with $J(0) = Y, J'(0) = UY $ satisfies $\|J(t) \| \leq C $ for $t \leq 0$.  The  subspace $E^u_{x,v}  \subset T_{x,v} T\M$ they generate corresponds to $TW^u_{x,v}$, where $W^u_{x,v}$  is the set of opposite vectors to vectors in $W^s_{x,-v}$. By definition, $S_{\dot \g_{x,}(0)} = - U_{\dot \g_{x,-v}(0)}$, so that we also have:
$$ B (x,v)\; := \;  - {\textrm {Tr }} S_{x,v} \; =\;  {\textrm {Tr }} U_{x,-v}.$$ 
We have ${\textrm {Ker }} S = {\textrm {Ker }} U$ and $Y \in {\textrm {Ker }}S $ if, and only if, the Jacobi field $J(t)$ with $J(0) = Y, J'(0) = 0 $ is bounded for all $t \in \R$. The rank of the geodesic $\g _{x,v} $  therefore is $\kappa = {\textrm { Dim Ker }} S$ and the geodesic $\g_{x,v} $ is of rank one only if Det($(U-S)|_{_{(\dot \g_{x,v}(0))^\perp )}} = 0$.

Recall that $SM $ is identified with the quotient of $\M \times \MM$ under the diagonal action of $G$. Clearly, for $g \in G$, $g(W^s_{x,v}) = W^s_{Dg (x,v)}$ so that the $W^s$ define a foliation $\W^s$ on $SM$. The leaves of the foliation $\W^s$ are quotient of $\M$, they are naturally endowed with the Riemannian metric induced from $\wt g$.

\

\

\subsection{Proof of Theorem \ref{main}}

We continue assuming that $(\M,\wt g)$ has nonpositive curvature. By the Rank Rigidity Theorem (see \cite {B3}), $(\M ,\wt g) $ is of the form $$( \M_0 \times \M_1 \times \cdots \times \M_j \times \M_{j+1} \times \cdots \times \M_k , \wt g)\footnote{With a clear convention for the cases  when Dim $\M_0 =0,$  $j =0 $ or   $k=j$.} ,$$
where $\wt g$ is the product metric $\wt g^2 = (\wt g_0)^2 + (\wt g_1)^2 + \cdots +  (\wt g_j)^2 +  (\wt g_{j+1})^2 +\cdots +  (\wt g_k)^2 $, $(\M_0 , \wt g_0) $ is Euclidean, $(\M_i, \wt g_i) $ is an irreducible symmetric space of rank at least two for $i = 1, \cdots, j$ and a rank-one manifold  for $i = j+1, \cdots , k.$ If the  $(\M_i, \wt g_i), i = j+1, \cdots k, $ are all symmetric spaces of rank one, then $(\M, \wt g)$ is a symmetric space. Moreover in that case, all inequalities in (\ref {ineq}) are equalities: this is the case for irreducible symmetric spaces
(all numbers are 0 for Euclidean space;  for the other spaces, $4\l_0 $ and $v^2$ are classically known to coincide  (\cite {O}) and we have:
$$ 4\l_0 (\M) \; = \; \sum _i 4 \l_0 (\M_i ), \quad v^2 (\M) \; = \; \sum _i v^2 (\M_i).$$
To prove Theorem \ref{main}, it suffices to prove that if $\ell ^2 = h $, all $\M_i $ in the decomposition are symmetric spaces. This is already true for $i = 0, 1, \cdots j.$ It remains to show that $(\M_i, \wt g_i )$ are symmetric spaces for $i = j+1, \cdots k$. Eberlein showed that each one of the spaces $(\M_i, \wt g_i) $ admits a cocompact discrete group of isometries (see \cite{Kn}, Theorem 3.3). This shows that  the linear drifts $\ell _i $ and the stochastic entropies $h_i $ exist for each one of the spaces $(\M_i, \wt g_i) $. Moreover, we clearly have
$$ \ell ^2 \; = \; \sum _i \ell _i^2, \quad h\; = \; \sum h_i .$$
Therefore Theorem \ref{main} follows from
\begin{theo} \label{rankone} Assume $(M,g)$ is a closed connected  rank one manifold of nonpositive curvature and that $\ell ^2 = h$. Then $(\M, \wt g)$ is a symmetric space. \end{theo}

\

A Hadamard manifold $\M$ is called asymptotically harmonic if the function $B (= \D_xb)$ is constant on $S\M$. Theorem \ref{rankone} directly follows from  two propositions:

\begin{prop} \label{as.harm} Assume $(M,g)$ is a closed connected  rank one manifold of nonpositive curvature and that $\ell ^2 = h$. Then $(\M, \wt g)$ is asymptotically harmonic. \end{prop}

\begin{prop}\label{final} [\cite{Z}, Theorem 1.1] Assume $(M,g)$ is a closed connected  rank one manifold of nonpositive curvature such  that   $(\M, \wt g)$ is asymptotically harmonic. Then, $(\M, \wt g)$ is
a symmetric space. \end{prop}

\

\section{Proof of Proposition \ref{as.harm}}

We consider the foliation $\W$ of subsection 2.2. Recall that the leaves are endowed with a natural Riemannian metric. We write $\D^\W $ for the associated Laplace  operator on functions which are of class $C^2$ along the leaves of $\W$. A probability measure $m$ on $SM$ is called harmonic if it satisfies, for any $C^2$ function $f$, we have:
$$ \int _{SM} \D^\W f  dm \; = \; 0 .$$
Let $M$ be a closed connected manifold such that $\ell ^2 = h$. In \cite{L3} it is shown that then, there exists a harmonic probability measure $m$ on $SM$ such that, at $m$-a.e. $(x,v)$, $B(x,v) = \ell $. Since $B$ is a continuous function, Proposition \ref{as.harm} follows from

\begin{theo} \label{harmo} Let $(M,g)$ be a closed connected rank one manifold of nonpositive curvature, $\W$ the stable foliation on $SM$ endowed with the natural metric as above. Then, there is only one harmonic probability  measure $m$ and the support of $m$ is the whole space $SM$.
\end{theo} 
\begin{proof}
Let $m$ be a $\W$ harmonic probability measure on $SM$. Then, there is a unique $G$-invariant measure $\wt m$ on $S\M$ which coincide with $m$ locally. Seen as a measure on $\M \times \MM$, we claim that $\wt m$ is given, for any $f$ continuous with compact support, by:
\begin{equation}\label{harm.meas} \int f(x, \xi ) d\wt m (x,\xi) \; = \; \frac {1}{{\textrm {Vol}} M} \int _{\M} \left( \int _{\MM} f(x,\xi ) d\nu _x (\xi) \right) dx, \end{equation}
where the family $x \mapsto \nu_x$ is a family of probability measures on $\MM$ such that, for all $\varphi $ continuous on $\MM$, $x\mapsto \int \varphi (\xi ) d\nu_x (\xi)  $ is a harmonic function on $\M$ and the measure $dx$ is  the Riemannian volume on $\M$. The claim follows from \cite{Ga}. For convenience,
let us reprove it: on the one hand, the measure $\wt m$ projects on $\M$ as a $G$-invariant measure satisfying $\int \D f dm  = 0 $. The projection of $\wt m $ on $\M$ is proportional to Volume, gives measure 1 to fundamental domains and  formula (\ref{harm.meas})  is the desintegration formula. On the other hand, if one projects $\wt m$ first on $\MM$, there is a probability measure $\nu $ on $\MM$ such that $$\int f(x,\xi) d\wt m (x,\xi) \;= \int_{\MM} \left( \int_{\M} f(x,\xi) dm_\xi (dx) \right) d\nu (\xi).$$
For $\nu $-a.e. $\xi $, the measure $m_\xi $ is a harmonic measure on $\M$; therefore, for $\nu $-a.e. $\xi $,  there is a positive harmonic function $k_\xi (x) $ such that $m_\xi = k_\xi (x) {\textrm{Vol.}}$ Comparing the two expressions for $\int f d\wt m$, we see that  the measure $\nu_x $ is given by $$ \nu _x \; = \; k_\xi (x) \nu$$ and $x \mapsto \int _{\MM} \varphi (\xi ) d\nu_x (\xi ) $ is indeed a harmonic function.

\

The $G$-invariance of $\wt m$ implies that, for all $g \in G$, $g_\ast \nu _x = \nu _{gx}$. In particular, the support of $\nu $ is $G$-invariant. By \cite{E1} (see \cite{B3}, page 48), the support of $\nu $ is the whole $\MM$ and therefore the support of $m$ is the whole $SM$. This result would be sufficient for proving Proposition \ref{as.harm}, but using discretization, we are going to  identify the measure $\nu _x$ on $ \M (\infty )$ as the hitting measure of the Brownian motion on $\M$ starting from $x$. This shows  Theorem \ref{harmo}.

\

Fix $x_0 \in \M$. The discretization procedure of Lyons and Sullivan (\cite{LS}) associates to the Brownian motion on $\M$ a probability measure $\mu $ on $G$ such that  $ \mu (g) >0 $ for all $g$ and that any bounded harmonic function $F$ on $\M$ satisfies $$ F(x_0) \; = \; \sum _{g \in G} F(gx_0 ) \mu (g).$$
Recall that for all $\varphi $ continuous on $\MM$, $x \mapsto \nu _x (\varphi ) $ is a harmonic function and that $\nu _{gx} = g_\ast \nu _x$. It follows that the measure $\nu _{x_0} $ is stationary for $\mu$, i.e. it satisfies:
$$ \nu _{x_0} \; = \; \sum _{g\in G} g_\ast \nu _{x_0} \mu (g).$$
Since the support of $\mu $ generates $G$ as a semigroup (actually, it is already the whole $G$), there is only one stationary probability measure on $\MM$ (see \cite {B3}, Theorem 4.11 page 58). We know one already: the hitting measure $m_{x_0} $ of the Brownian motion on $\M$ starting from $x_0$. This shows that $\nu _{x_0} = m_{x_0}$. Since $x_0$ was arbitrary in the above reasoning, we have $\nu _x = m_x $ for all $x \in \M$ and the measure $\wt m$ is given by:
$$\int f(x, \xi ) d\wt m (x,\xi) \; = \; \frac {1}{{\textrm {Vol}} M} \int _{\M} \left( \int _{\MM} f(x,\xi ) dm_x (\xi) \right) dx. $$
\end{proof}

\

 {\bf {Acknowledgements}} I am   very grateful to Gerhard Knieper for his interest and his comments, in particular for having attracted my attention to \cite{Z}. I also acknowledge partial support of NSF grant DMS-0811127.

\


\small\begin{thebibliography}{99}
\bibitem[{\bf B1}]{B1} W. Ballmann, Nonpositively curved manifolds of higher rank, {\em Ann. Math}, {\bf 122} (1985), 597--609.

\bibitem[{\bf B2}]{B2} W. Ballmann, On the Dirichlet problem at infinity for manifolds of nonpositive curvature, {\em Forum Mathematicum}, {\bf 1} (1989), 201--213.

\bibitem[{\bf B3}]{B3} W. Ballmann, Lectures on spaces of nonpositive curvature, {\em DMV Seminar}, {\bf 25} (1995).

\bibitem[{\bf  BCG}]{BCG} G. Besson, G. Courtois and S. Gallot, Entropies et rigidit\'es des espaces localement sym\'etriques de courbure strictement n\'egative, {\em Geom. Func. Anal.}   {\bf 5} (1995), 731--799.

\bibitem[{\bf BFL}]{BFL} Y. Benoist, P. Foulon and F. Labourie,  Flots d'Anosov \`a distributions stables et instables diff\'erentiables, {\em  J. Amer. Math. Soc. }{\bf 5} (1992), 33--74.

\bibitem[{\bf BS}]{BS} K. Burns and R. Spatzier, Manifolds of nonpositive curvature and their buildings, {\em Publications math. IHES}, {\bf 65} (1987),  35--59.

\bibitem[{\bf E}]{E1} P. Eberlein, Geodesic flows on negatively curved manifolds, II, {\em Transactions Amer. math. Soc.}, {\bf 178} (1973), 57--82.

\bibitem[{\bf  FL}]{FL}   P. Foulon and F. Labourie,  Sur les vari\'et\'es compactes asymptotiquement harmoniques, {\em  Invent. Math. }{\bf 109} (1992), 97--111.
 
\bibitem[{\bf Ga}]{Ga} L. Garnett,   Foliations, the ergodic theorem and Brownian motion, {\em J. Funct. Anal. } {\bf 51} (1983), 285--311.


\bibitem[\bf Gu] {Gu}Y. Guivarc'h, Sur la loi des grands nombres et le rayon spectral
d'une marche al\'{e}atoire, \emph{Ast\'{e}risque,} {\bf 74} (1980) 47--98.

\bibitem[{\bf HI}]{HI}   E. Heintze and H.-C. Im Hof,  Geometry of horospheres,  {\em  J. Diff. Geom.} {\bf 12} (1977), 481--491.

 \bibitem [\bf K]{K}V. A. Kaimanovich, Brownian motion and harmonic
functions on covering manifolds. An entropic approach, {\em Soviet Math.
Dokl.} {\bf 33} (1986) 812--816.


\bibitem[{\bf Ka}]{Ka} A. Katok, Four applications of conformal equivalence to geometry and dynamics, {\em Ergod. Th. \& Dynam. Sys.}, {\bf $8^\ast $} (1988), 139--152.



\bibitem[{\bf Kn}]{Kn} G. Knieper, On the asymptotic geometry of non-positively curved manifolds, {\em GAFA}, {\bf 7} (1997), 755--782.

\bibitem[{\bf L1}]{L1} F. Ledrappier, Harmonic measures and Bowen-Margulis measures, {\em Israel J. Math.} {\bf71} (1990), 275--287.

\bibitem[{\bf L2}]{L2} F. Ledrappier, Applications of dynamics to compact manifolds of negative
curvature,
{\it in}  {\sl Proceedings of the ICM Z\"urich 1994}, Birkh\"auser (1995),
1195-1202.


\bibitem[{\bf L3}]{L3} F. Ledrappier, Linear drift and entropy for regular covers, {\em GAFA}, {\bf 20} (2010), 710--725.


\bibitem[{\bf LS}]{LS} T. Lyons and D. Sullivan, Function theory, random paths and covering spaces, {\em J. Differential Geometry}, {\bf 19} (1984), 299--323.

\bibitem[{\bf O}]{O} M.A. Olshanetsky, Martin boundary for the Laplace-Beltrami operator on a Riemannian symmetric space of non-positive curvature, {\em Uspehi Mat. Nauk.,} {\bf 24:6} (1969), 189-190.

\bibitem[{\bf Z}]{Z} A. M. Zimmer, Asymptotically harmonic manifolds without focal points, {\em preprint}, (http://arxiv.org/abs/1109.2481),

\end{thebibliography}
\end{document}